\documentclass[10pt]{article}

\usepackage{theorem,amssymb,amsmath}

\usepackage{graphicx}

\usepackage{color}                    
\definecolor{red}{rgb}{1,0,.2}        
\definecolor{cjp}{rgb}{.1,.7,.2}        
\definecolor{fmdc}{rgb}{1,0,.8}        

\newcommand{\T}{x_\beta}

\newtheorem{theorem}{Theorem}

\newtheorem{lemma}[theorem]{Lemma}

\newtheorem{proposition}[theorem]{Proposition}
\theorembodyfont{\rm}

\advance \topmargin by -\headheight
\advance \topmargin by -\headsep
\evensidemargin \oddsidemargin
\begin{document}

\bigskip

\begin{centering}

{\bf \huge The sum of digits function of the base \\[.5cm] phi expansion of the natural numbers}

\vspace*{1cm}

{\bf \Large F.~Michel Dekking}

\bigskip

\bigskip

{\large \bf \it Dedicated to the memory of Christian Mauduit}

\bigskip

\bigskip

{ F.M.Dekking@TUDelft.nl}

\medskip

{ DIAM,  Delft University of Technology, Faculty EEMCS,\\ P.O.~Box 5031, 2600 GA Delft, The Netherlands.}

\bigskip

\begin{abstract}
{\small In the base phi expansion any natural number is written uniquely as a sum of powers of the golden mean with digits 0 and 1, where one requires that the product of two consecutive digits is always 0. In this paper we show that the sum of digits function modulo 2 of these expansions is a morphic sequence. In particular we prove that  --- like for the Thue-Morse sequence --- the frequency of 0's and 1's in this sequence is equal to 1/2.}
\end{abstract}

\medskip



\date{\today}

\end{centering}

\bigskip

\section{Introduction}

  Base phi representations were introduced by George Bergman in 1957 (\cite{Bergman}). Base phi representations are also known  as beta-expansions of the natural numbers, with $\beta=(1+\sqrt{5})/2=:\varphi$, the golden mean.\\
  A natural number $N$ is written in base phi if $N$ has the form
  $$N= \sum_{i=-\infty}^{\infty} d_i \varphi^i,\vspace*{-.0cm}$$
  with digits $d_i=0$ or 1, and where $d_id_{i+1} = 11$ is not allowed.
  We  write these expansions as
  $$\beta(N) = d_{L}d_{L-1}\dots d_1d_0\cdot d_{-1}d_{-2} \dots d_{R+1}d_R.$$

  \noindent Ignoring leading and trailing 0's, the base phi representation of a number $N$ is unique (\cite{Bergman}).

  Let for $N\ge 0$
  $$s_\beta(N):=\sum_{k=L}^{k=R} d_k(N)$$
  be the sum of digits function of the base phi expansions. We have
  $$(s_\beta(N)) = 0, 1, 2, 2, 3, 3, 3, 2, 3, 4, 4, 5, 4, 4, 4, 5, 4, 4, 2, 3, 4, 4, 5, 5, 5, 4, 5, 6, 6, 7, 5, 5, 5, 6, 5, 5, 4, 5, 6, 6, 7,\dots$$
  In \cite{Cooper} asymptotic expressions as $x\rightarrow\infty$ for $\sum_{N<x}s_\beta(N)$ where obtained.

  \medskip

  \noindent In this paper we  study the base phi analogue of the Thue-Morse sequence (where the base equals 2), i.e., the sequence
  $$(\T(N)):= (s_\beta(N) \!\!\! \mod 2) \,= \, 0,1,0,0,1,1,1,0,1,0,0,1,0,0,0,1,0,0,0,1,0,0,1,1,1,0,1,0,0,\dots$$

  \medskip

  Recall that a morphism is a map from the set of infinite words over an alphabet to itself, respecting the concatenation operation.
The Thue Morse sequence is the fixed point starting with 0 of the morphism $0\rightarrow 01,\;1\rightarrow 10$.

\bigskip

\noindent {\bf Theorem}\label{th:main} {\it\; $\T$ is a morphic sequence, i.e., the letter-to-letter image of the fixed point of a morphism.}

\bigskip

This theorem permits to answer a number of natural questions one may ask about $\T$, as for example: will a word 00000 ever occur? What are the frequencies of 0 and 1?

\section{Properties of the base phi representation}\label{sec:morph}

The Lucas numbers $(L_n)=(2, 1, 3, 4, 7, 11, 18, 29, 47, 76,123, 199, 322,\dots)$ are defined by
$$  L_0 = 2,\quad L_1 = 1,\quad L_n = L_{n-1} + L_{n-2}\quad {\rm for \:}n\ge 2.$$
The Lucas numbers have a particularly simple base phi representation: from  the well-known formula
$L_{2n}=\varphi^{2n}+\varphi^{-2n}$, and the recursion $L_{2n+1}=L_{2n}+L_{2n-1}$, we have for all $n\ge 1$
$$ \beta(L_{2n}) = 10^{2n}\cdot0^{2n-1}1,\quad \beta(L_{2n+1}) = 1(01)^n\cdot(01)^n.$$

\medskip

\noindent The properties of base phi expansion of the natural numbers can be read of from the {\it Lucas intervals}:
 $$\Lambda_{2n}:=[L_{2n},L_{2n+1}], \quad \Lambda_{2n+1}:=[L_{2n+1}+1, L_{2n+2}-1].$$
 When we add $\Lambda_0:=[0,1]$, these intervals partition the natural numbers as $n=0,1,2\dots$. The partition elements correspond to the lengths of the expansions: if 
 $\beta(N) = d_{L}d_{L-1}\dots d_1d_0\cdot d_{-1}d_{-2} \dots d_{R+1}d_R,$
then the left most index $L=L(N)$ and the right most index $R=R(N)$ satisfy
$$L(N)=2n\!+1, \,R(N)=-2n \;{\rm iff}\; N\in \Lambda_{2n}, \quad L(N)=2n\!+2= -R(N) \;{\rm iff}\; N\in \Lambda_{2n+1}.$$
This is not hard to see from the simple expressions we have for the $\beta$-expansions of the Lucas numbers, see also Theorem 1 in \cite{Grabner94}.

\medskip

\noindent Since $\beta(L_{2n})$ consists of only 0's between the exterior 1's, the following lemma is obvious.

\begin{lemma}\label{lem:even} {\rm \bf (\cite{Dekking-FQ})} For all $n\ge 1$ and $k=0,\dots,L_{2n-1}$
one has $ \beta(L_{2n}+k) =  \beta(L_{2n})+ \beta(k).$
\end{lemma}

  This gives recursive relations for the expansions in the Lucas interval $\Lambda_{2n}$. To obtain recursive relations for the interval $\Lambda_{2n+1}$, this interval has to be divided into three subintervals. These three intervals are
$$I_n:=[L_{2n+1}+1,\, L_{2n+1}+L_{2n-2}-1],\: J_n:=[L_{2n+1}+L_{2n-2},\, L_{2n+1}+L_{2n-1}], \: K_n:=[L_{2n+1}+L_{2n-1}+1,\, L_{2n+2}-1].$$


\medskip

To formulate the following lemma, it is notationally convenient to extend the semigroup of words to the free group of words. For example, one has $110^{-1}01^{-1}00=100$.

\begin{lemma}\label{lem:odd}{\rm \bf (\cite{San-San}, \cite{Dekking-FQ})}  For all $n\ge 2$ and $k=1,\dots,L_{2n-2}-1$
\begin{align*}
I_n:&\quad \beta(L_{2n+1}+k) = 1000(10)^{-1}\beta(L_{2n-1}+k)(01)^{-1}1001,\\ K_n:&\quad\beta(L_{2n+1}+L_{2n-1}+k)=1010(10)^{-1}\beta(L_{2n-1}+k)(01)^{-1}0001=10\beta(L_{2n-1}+k)(01)^{-1}0001.
\end{align*}
Moreover, for all $n\ge 2$ and $k=0,\dots,L_{2n-3}$
$$\hspace*{-3.8cm}J_n:\quad\beta(L_{2n+1}+L_{2n-2}+k) = 10010(10)^{-1}\beta(L_{2n-2}+k)(01)^{-1}001001.$$
\end{lemma}

\section{The sequence $\T$ is morphic}\label{sec:proof}

If $V=[K,K\!+1,\dots,L]$ is an interval of natural numbers, then we write
$$\T(V):=[\T(K),\T(K\!+1),\dots,\T(L)]$$
 for the  consecutive sums of digits modulo 2 of these numbers.

\noindent Since $\T(L_{2n})=0$ and $\T(0) = 0$, Lemma \ref{lem:even} implies directly the following lemma.     

\begin{lemma}\label{lem:EVEN} {\rm(LEMMA EVEN)} For all $n\ge 1$ one has $\T(\Lambda_{2n}) =  \T([0,L_{2n-1}]).$
\end{lemma}

\noindent The mirror morphism on $\{0,1\}$ is defined by $\overline{0}=1, \overline{1}=0$.

\medskip

\noindent We obtain from Lemma \ref{lem:odd} with $\T(I_n)=\T(K_n)=\overline{ \T(\Lambda_{2n-1})}$, and   $\T(J_n)= \T(\Lambda_{2n-2})$:

\begin{lemma}\label{lem:ODD} {\rm(LEMMA ODD)} For all $n\ge 1$ one has
$\T(\Lambda_{2n+1}) = \overline{ \T(\Lambda_{2n-1})} \T(\Lambda_{2n-2})\overline{\T(\Lambda_{2n-1})}.$
\end{lemma}

\medskip

\noindent We illustrate the base phi expansions with the following table.

\medskip

 \begin{tabular}{|r|c|c|c|}
   \hline
  \;$N^{\phantom{|}}$ & $\beta(N)$ &\; $\T(N)$ & {\small Lucas interval}\\[.0cm]
   \hline
   0\;& \;\!\!\!$0$          & $0$  & $\Lambda_0$ \\
   1\;& \;\!\!\!$1$          & $1$  & $\Lambda_0$ \\
   \hline
   2\; & \;\:\,\,$10\cdot01$  & $0$  & $\Lambda_1$\\
   \hline
   3\; & \;$100\cdot01$     & $0$    & $\Lambda_2$\\
   4\; & \;$101\cdot01$     & $1$ & $\Lambda_2$\\
   \hline
   5\; & \;\:\,$1000\cdot1001$  & $1$ & $\Lambda_3$\\
   6\; & \;\:\,$1010\cdot0001$  & $1$ & $\Lambda_3$\\
    \hline
   7\; & \;$10000\cdot0001$ & $0$ & $\Lambda_4$\\
   8\; & \;$10001\cdot 0001$  & $1$ & $\Lambda_4$\\
   9\; & \;$10010\cdot0101$     & $0$ & $\Lambda_4$\\
   10\; & \;$10100\cdot0101$     & $0$ & $\Lambda_4$\\
   11\; & \;$10101\cdot0101$     & $1$ & $\Lambda_4$\\
    \hline
   12\; & \;\,\,$100000\cdot101001$  & $0$ & $\Lambda_5$\\
   \hline
 \end{tabular}  \quad

 \bigskip

\noindent Let $\tau$ be the morphism on the alphabet $A:=\{1,\dots,8\}$ defined by\\[-.6cm]
\begin{align*} 
 \tau(1) = & \, 12,  \qquad  \tau(2) = \, 312,  \qquad  \tau(3) = \, 47,  \qquad \tau(4) = \, 8312, \\[-.1cm]
 \tau(5) = & \,56,   \qquad  \tau(6) = \,756,   \qquad  \tau(7) = \, 83,  \qquad \tau(8) = \, 4756.
\end{align*} 
Define the mirroring morphism $\mu$ on $A$ by
$$\mu:\; 1\rightarrow 5,\; 2\rightarrow 6,\; 3\rightarrow 7,\; 4\rightarrow 8,\; 5\rightarrow 1,\; 6\rightarrow 2,\; 7\rightarrow 3,\; 8\rightarrow 4.$$
Then $\tau$ is mirror invariant:   $\tau\mu=\mu\tau$.

\begin{theorem}\label{th:fix} Let $\T$ be the sum of digits function of the base phi expansions of the natural numbers.

\noindent Let $\lambda:A^*\rightarrow\{0,1\}$ be the letter-to-letter morphism given by
$$\lambda(1)=\lambda(3)=\lambda(6)=\lambda(8)=0,\;{\rm  and}\;
\lambda(2)=\lambda(4)=\lambda(5)=\lambda(7)=1.$$
Then $\T=\lambda(t)$, where $t=1231247123\dots$ is the fixed point of $\tau$ starting with 1.
\end{theorem}

\medskip

\noindent Theorem \ref{th:fix} is a direct consequence of the following result.  Note that $\overline{\lambda\tau}=\lambda\tau\mu$.

\begin{proposition} For $n=1,2\dots$ one has $\T(\Lambda_{2n})=\lambda(\tau^n(1))$, and $\T(\Lambda_{2n+1})=\lambda(\tau^n(3))$.
\end{proposition}

\noindent {\it Proof:} By induction.
For $n=1$ one has $\T(\Lambda_{2})=01=\lambda(12)=\lambda(\tau(1))$, and $\T(\Lambda_{3})=11=\lambda(47)=\lambda(\tau(3))$.

\noindent From Lemma \ref{lem:EVEN} and the induction hypothesis we have\\[-.6cm]
\begin{align*}
\T(\Lambda_{2n+2}) &=  \T([0,L_{2n-1}])\T([L_{2n-1}+1,L_{2n}-1])\T([L_{2n},L_{2n+2}])\\
&=\lambda(\tau^n(1))\lambda(\tau^{n-1}(3))\lambda(\tau^n(1))\\
&= \lambda(\tau^{n-1}(12312))= \lambda(\tau^{n+1}(1)).
\end{align*}
\noindent From Lemma \ref{lem:ODD} and the induction hypothesis we have\\[-.6cm]
\begin{align*}
\T(\Lambda_{2n+3})&=\overline{\T(\Lambda_{2n+1})}\T(\Lambda_{2n})\overline{\T(\Lambda_{2n+1})}\\
 &=  \overline{\lambda(\tau^n(3)}\lambda(\tau^n(1))\overline{(\lambda(\tau^n(3))}\\
&=\lambda(\tau^n(7))\lambda(\tau^{n}(1))\lambda(\tau^n(7))\\
&=\lambda(\tau^n(717))=\lambda(\tau^n(47))=\lambda(\tau^{n+1}(3)). \hspace*{8cm} \Box
\end{align*}

\medskip

Since $\tau$ is mirror invariant, the letters $a$ and $\mu(a)$ have the same frequency for $a\in A$.
As $\overline{\lambda}=\lambda\mu$, this implies the following.

\begin{proposition} The letters 0 and 1 have frequency $\tfrac12$ in $\T$.
\end{proposition}

It is well-known that the words of length 2 in the Thue-Morse sequence have frequencies  $\tfrac16$ for 00 and 11, and $\tfrac13$ for 01 and 10.
Here is the corresponding result for the golden mean sum of digits function.

\begin{proposition} In $\T$ the words 00 and 11 have frequency $\tfrac1{10}\sqrt{5}$, and the words 01 and 10 have frequency $\tfrac12-\tfrac1{10}\sqrt{5}$.
\end{proposition}

\noindent {\it Proof:} As in \cite{Queff} we compute the frequencies $\nu[ab]$ of the words $ab$ of length 2 occurring in  the fixed point $t$ of the morphism $\tau$ by using the 2-block substitution $\tau^{[2]}$.
The words of length 2 occurring in the fixed point $t$ of the morphism $\tau$ are\\[-.6cm]
$$12,\: 23,\: 24,\: 28,\: 31,\: 35,\: 47,\: 56,\: 64,\: 67,\: 68,\: 71,\: 75,\: 83.$$
When we code the 14 words of length 2 by $\ell_1,\dots,\ell_{14}$, in the order given above, then $\tau^{[2]}$ is given for the letters $\ell_1,\dots\ell_7$ by
$$\tau^{[2]}:\quad \ell_1\rightarrow \ell_1\ell_2,\: \ell_2 \rightarrow \ell_5\ell_{13},\: \ell_3\rightarrow \ell_5\ell_{14},\: \ell_4\rightarrow \ell_5\ell_{13},\: \ell_5\rightarrow \ell_7\ell_{12},\: \ell_6\rightarrow \ell_7\ell_{13},\: \ell_7\rightarrow \ell_{14}\ell_5\ell_{14}.$$
The $\tau^{[2]}$-images of $\ell_8,\dots,\ell_{14}$ follow from this by mirror-symmetry. The first 7 components of the normalized eigenvector of the incidence matrix of the morphism $\tau^{[2]}$ are given by

$$\Big[\tfrac14-\tfrac1{20}\sqrt{5},\: \tfrac12-\tfrac15\sqrt{5},\: \tfrac3{20}-\tfrac1{20}\sqrt{5},\: \tfrac15\sqrt{5}-\tfrac25,\: \tfrac1{10},\: \tfrac3{20}-\tfrac1{20}\sqrt{5},\:\: \tfrac3{20}\sqrt{5}-\tfrac14\Big].$$

\noindent This means that, e.g., $\nu[12]=\tfrac14-\tfrac1{20}\sqrt{5}$, and $\nu[31]=\tfrac1{10}$.

\noindent The frequency of 00 equals  $\mu[00]=\nu[13]+\nu[68]+\nu[83]=\tfrac1{10}\sqrt{5}$. \hfill $\Box$

\bigskip

\noindent {\bf Final remark} Christian Mauduit with Michael  Drmota and  Jo\"el  Rivat proved that the Thue-Morse sequence is normal along squares (see \cite{Drm-Mau-Riv}).
Conjecture: this also holds for the sum of digits function modulo 2 of the basis phi expansion of the natural numbers, i.e., for $(\T(n^2))$.

\section*{Acknowledgement}
I am grateful to the organizers of the CIRM conference  ``Prime Numbers, Determinism and Pseudorandomness" for creating an excellent environment for remembering Christian Mauduit and his work.

\noindent I also thank Peter Grabner for providing a relevant reference.

\end{document}